\documentclass[12pt]{article}

\newcommand{\ncm}{\newcommand}
\ncm{\aut}{auto\-mor\-phi\-sm} \ncm{\Inn}{\mbox{\rm Inn($A$)}} 
\ncm{\Ap}{\mbox{$\overline{\rm Inn}(A)$}} \ncm{\Ext}{\mbox{\rm 
Ext}} \ncm{\Ex}{\mbox{\rm Ex}} \ncm{\OExt}{\mbox{\rm OrderExt}} 
\ncm{\AI}{\mbox{\rm AInn($A$)}} \ncm{\HI}{\mbox{\rm HInn($A$)}} 
\ncm{\Aut}{\mbox{\rm Aut}} \ncm{\Mal}{\mbox{$M_{\alpha}$}} 
\ncm{\Aff}{\mbox{${\rm Aff}$}} \ncm{\id}{\mbox{\rm id}} 
\ncm{\Ker}{\mbox{\rm Ker}} \ncm{\BE}{\begin{eqnarray*}} 
\ncm{\EE}{\end{eqnarray*}} \ncm{\lra}{\mbox{$\longrightarrow$}} 
\ncm{\Hom}{\mbox{\rm Hom}} \ncm{\calU}{{\cal U}} \ncm{\el}{\ell} 
\ncm{\ad}{\mbox{\rm ad}} \ncm{\Alg}{\mbox{\rm Alg}} 
\ncm{\Conv}{\mbox{\rm Conv}} \ncm{\D}{{\cal D}} 

\ncm{\cstar}{$C^{*}$-algebra} \ncm{\cstars}{$C^{*}$-algebras} 
\ncm{\ra}{\mbox{$\rightarrow$}} \ncm{\la}{\mbox{$\leftarrow$}} 
\ncm{\hra}{\hookrightarrow} \ncm{\da}{\mbox{$\downarrow$}} 
\ncm{\se}{\mbox{$\searrow$}} \ncm{\al}{\mbox{$\alpha $}} 
\ncm{\del}{\mbox{$\delta$}} \ncm{\supp}{\mbox{\rm supp}} 
\ncm{\Ad}{\mbox{\rm Ad}} \ncm{\CAR}{\mbox{$M_{2^{\infty}}$}} 
\ncm{\ep}{\mbox{$\epsilon > 0$}} \ncm{\mod}{\mbox{\rm mod}} 
\ncm{\Sp}{\mbox{\rm Sp}} \ncm{\ol}{\overline} 
\ncm{\Mninf}{\mbox{$M_{n^{\infty}}$}} \ncm{\MR}{M. R\o{}rdam} 
\ncm{\Range}{\mbox{\rm Range}} 
\ncm{\vo}{}
\ncm{\ch}{}
\ncm{\CMP}{Comm. Math. Phys.} \ncm{\add}{} 
\ncm{\tilsig}{\tilde{\sigma}} \ncm{\dist}{{\rm 
dist}}\ncm{\eps}{\epsilon} 

\newtheorem{theo}{Theorem}[section]
\newtheorem{cor}[theo]{Corollary}
\newtheorem{lem}[theo]{Lemma}
\newtheorem{prop}[theo]{Proposition}
\newtheorem{remark}[theo]{Remark}
\newtheorem{definition}[theo]{Definition}
\newtheorem{example}[theo]{Example}

\newenvironment{rem}{\begin{remark} \rm}{\end{remark}}
\newenvironment{pf}{{\it Proof.}}{\vspace{3mm}}

\ncm{\R}{\mbox{\bf R}} \ncm{\Z}{\mbox{\bf Z}} \ncm{\T}{\mbox{\bf 
T}} \ncm{\TT}{\T$^{2}$} \ncm{\N}{\mbox{\bf N}} \ncm{\C}{\mbox{\bf 
C}} 

\title{UHF flows and the flip automorphism}
\bigskip

\author{A. Kishimoto\\
  \small Department of Mathematics, Hokkaido University,
   Sapporo 060, Japan}
\date{\small November 2000}
\begin{document}
\maketitle

\section{Introduction}
By a {\em flow} on a unital \cstar\ $A$ we mean a strongly 
continuous one-parameter automorphism group. The infinitesimal 
generator $\delta_{\alpha}$ of a flow $\alpha$ is a closed 
derivation in $A$, by which we mean that $\delta_{\alpha}$ is a 
closed linear operator which is defined on a dense *-algebra 
$D(\delta_{\alpha})$ and satisfies that 
$\delta_{\alpha}(x)^*=\delta_{\alpha}(x^*)$ and 
$\delta_{\alpha}(xy)=\delta_{\alpha}(x)y+x\delta_{\alpha}(y)$. See 
\cite{Br,BR,Sak} for the general theory of derivations. 

If $\alpha$ is a flow and $u=(u_t)$ is a one-parameter family of 
unitaries of $A$ such that $t\mapsto u_t$ is continuous and 
$u_s\alpha_s(u_t)=u_{s+t}$ for all $s,t\in\R$, we say that $u$ is 
an $\alpha$-cocycle.  Then $\Ad\,u\alpha:t\mapsto \Ad\,u\alpha_t$ 
is a flow and is called a cocycle perturbation of $\alpha$. If $u$ 
is differentiable with 
 $$h=-i\frac{du_t}{dt}|_{t=0},
 $$
then $\Ad\,u\alpha$ is an inner perturbation of $\alpha$ in the 
sense that $\delta_{{\rm Ad}\,u\alpha}=\delta_{\alpha}+\ad\,ih$, 
where $\ad\,ih$ is the inner derivation defined by 
$\ad\,ih(x)=i(hx-xh)$ for $x\in A$. We write $\Ad\,u\alpha$ also 
as $\alpha^{(h)}$ in this case. 

If $v$ is a unitary which is not in the domain 
$D(\delta_{\alpha})$, then the flow 
 $$\Ad\,v\,\alpha\,\Ad\,v^*:t\mapsto\Ad\,v\alpha_t(v^*)\,\alpha_t,
 $$ 
which is conjugate to $\alpha$, is a cocycle perturbation of 
$\alpha$ but not an inner perturbation. In general a cocycle 
perturbation is conjugate to an inner perturbation since any 
$\alpha$-cocycle $u$ is cohomologous to a differentiable one $w$, 
i.e., $u_t=vw_t\alpha_t(v^*)$ for some unitary $v$ (see 
\cite{K5}). 

By an inner flow we mean a flow given by $t\mapsto 
\Ad\,e^{iht}=e^{t\,{\rm ad}ih}$ for some self-adjoint element $h$ 
in $A$. We say that $\alpha$ is approximately inner if there is a 
sequence $(h_n)$ in $A_{sa}$ such that $\alpha_t=\lim 
\Ad\,e^{ih_nt}$, i.e., $\alpha_t(x)=\lim\Ad\,e^{ih_nt}(x)$ for 
every $t\in\R$ and $x\in A$, or equivalently, uniformly continuous 
in $t$ on every compact subset of $\R$ and every $x\in A$. 

When $A$ is an AF \cstar, we call $\alpha$ an {\em AF flow} if 
there is an increasing sequence $(A_n)$ of finite-dimensional 
$C^*$-subalgebras of $A$ such that $1_A\in A_1$, $A=\ol{\cup_n 
A_n}$ and $\alpha_t(A_n)=A_n$ for all $n$ and $t$. In this case 
there is an $h_n$ in $(A_n)_{sa}$ such that 
$\alpha_t|A_n=\Ad\,e^{ih_nt}|A_n$ and hence $\alpha$ is 
approximately inner. (The term AF flow is coined in \cite{BK2} 
while the adjective {\em locally representable} is used in 
\cite{K5} and elsewhere to refer to the same object.)  

When $A$ is a UHF \cstar, we call $\alpha$ a {\em UHF flow} if 
there is an increasing sequence $(A_n)$ of full matrix 
$C^*$-subalgebras of $A$ such that $1_A\in A_1$, 
$A=\ol{\cup_nA_n}$, and $\alpha_t(A_n)=A_n$ for all $n$ and $t$.  
If we set $B_n=A_n\cap A_{n-1}'$ with $A_0=0$, then $A\cong 
\bigotimes_{1}^{\infty}B_n$ and $\alpha_t(B_n)=B_n$ for all $n$. 
Namely $\alpha$ is of infinite tensor product type. Note that UHF 
flows are AF flows and there are AF flows on UHF \cstars\ which 
are not UHF flows. This follows because there are AF flows which 
have more than one KMS states for some temperature while UHF flows 
always have a unique KMS state for any temperature (see e.g., 
\cite{K5}). 

We call $\alpha$ a {\em compact flow} if the closure of 
$\alpha_{\bf R}$ in $\Aut(A)$, the automorphism group of $A$, is 
compact, or equivalently, if there is an increasing sequence 
$(V_n)$ of finite-dimensional subspaces of $A$ such that 
$A=\ol{\cup_nV_n}$ and $\alpha_t(V_n)=V_n$ for all $n$ and $t$. A 
periodic flow is compact but there are more. If $\alpha$ is an AF 
flow, then it is compact. 

In \cite{K6} we refer to locally inner (or locally representable) 
flows as a generalization of AF flows. That is, a flow  $\alpha$ 
is {\em locally inner} if there is an increasing sequence $(A_n)$ 
of (arbitrary) $C^*$-subalgebras of $A$ such that $A=\ol{\cup_n 
A_n}$ and $\alpha$ leaves $A_n$ invariant and restricts to an 
inner flow on $A_n$. Apparently AF flows are locally inner and 
locally inner flows are approximately inner.  But there seems to 
be no obvious relation between compact flows and locally inner 
flows.    

When $A$ is a UHF \cstar, we have the following implications for 
flows: 
 \begin{quote} UHF flows $\Rightarrow$ AF flows 
  $\left\{\begin{array}{ll}\Rightarrow& Locally\ inner\ flows\\
  \Rightarrow& C ompact\ flows\end{array}\right.$  
 \end{quote}
where the reverse implications are false (for the latter  note 
that the examples of non-AF flows constructed in \cite{K6} are 
locally inner and compact). 
  
There is certainly a non-compact flow (e.g., a flow with some 
asymptotic abelianess), but we do not seem to know if there is a 
flow which is not even a cocycle perturbation of a compact flow. 
If there is a flow which is not approximately inner (against the 
Powers-Sakai Conjecture \cite{Sak}), it is likely that we have 
such a one among the compact (or even periodic) flows since there 
is such a one for some simple AF \cstar, which is obtained from 
\cite{K7} and Lin's classification theorem for simple \cstars\ of 
tracial topological rank zero. Hence there may be no inclusion 
relation between the compact flows and the approximately inner 
flows. 

In this note  we will briefly discuss compact flows and then look 
into UHF flows, a simplest kind of flows! 

We will show in \ref{D2} that if $\alpha$ is an approximately 
inner compact flow on $A$, then the domain $D(\delta_{\alpha})$ of 
the generator $\delta_{\alpha}$ contains a maximal abelian 
$C^*$-subalgebra (masa for short) of $A$, a property which 
obviously holds for locally inner flows. We here refer to a result 
that the AF flows are characterized by the property that the 
domain contains a canonical AF masa (see \cite{K6} for details). 

When $A$ is a UHF \cstar\ such that $A\otimes A\cong A$, we call a 
UHF flow $\gamma$ on $A$ {\em universal} if $\gamma\otimes\alpha$ 
is cocycle conjugate to $\gamma$ for any UHF flow $\alpha$ on $A$ 
(i.e., there is an isomorphism $\varphi$ of $A\otimes A$ onto $A$ 
such that $\gamma\otimes\alpha$ is a cocycle perturbation of 
$\varphi^{-1}\gamma\varphi$). There exist universal UHF flows on 
$A$ (as shown in \ref{A1}) and they are mutually cocycle 
conjugate, or even {\em almost conjugate} (see \ref{B4}). In the 
case of the CAR algebras we will construct a universal UHF flow in 
a simple way extending a result given in \cite{K5}; see \ref{C1}. 

When $\alpha$ is a flow on $A$ and $\sigma$ is the flip 
automorphism of $A\otimes A$, i.e., $\sigma(x\otimes y)=y\otimes 
x,\ x,y\in A$, we say that $\sigma$ is {\em $\alpha$-invariantly 
approximately inner} if there is a sequence $(u_n)$ of unitaries 
in $A\otimes A$ such that $\sigma=\lim\Ad\,u_n$ and 
$\|(\alpha_t\otimes\alpha_t)(u_n)-u_n\|\ra0$ for each $t\in\R$, or 
equivalently uniformly in $t$ on every compact subset of $\R$. We 
shall prove for an approximately inner flow $\alpha$ that 
$\alpha\otimes\gamma$ is a cocycle perturbation of a (universal) 
UHF flow if and only if $\sigma$ is $\alpha$-invariantly 
approximately inner, where $\gamma$ is a universal UHF flow; see 
\ref{B5}. We do not really have any application of this result 
(except for a result on quasi-UHF flows; see \ref{B6}) but this is 
a first attempt to characterize UHF flows; see \cite{K6,BK2} for 
some results on AF flows.  We shall also note that if the flip is 
$\alpha$-smoothly approximately inner, i.e., if we replace the 
condition $\|(\alpha_t\otimes\alpha_t)(u_n)-u_n\|\ra0$ by that 
$(\alpha_t\otimes\alpha_t(u_n))$ is equi-continuous in the above 
definition, then $\alpha$ has a unique KMS state for each inverse 
temperature, see \ref{B3}. This follows from \cite{FVV} and 
applies to the one-dimensional quantum lattice systems, thus 
providing another proof of the well-known uniqueness result.

\section{Compact flows}
\begin{prop}\label{D1}
Let $A$ be a separable \cstar\ and $\alpha$ a flow on $A$. Then 
the following conditions are equivalent: 
 \begin{enumerate}
 \item The closure of $\{\alpha_t\,|\,t\in\R\}$ in $\Aut(A)$ is compact.
 \item For each $x\in A$ the closure of $\{\alpha_t(x)\,|\,t\in\R\}$ is
 compact.
 \item For each $x$ of a dense subset of $A$ the closure of 
 $\{\alpha_t(x)\,|\,t\in\R\}$ is compact.
 \item The linear span of $A_p=\{x\in A\ |\ 
 \alpha_t(x)=e^{ipt}x\}$ for all $p\in\R$ is dense in $A$.
 \item There is an increasing sequence $(V_n)$ of finite-dimensional 
 subspaces of $A$ such that $A=\ol{\cup_nV_n}$ and $\alpha_t(V_n)=V_n$ for 
 all $n$ and $t$.
 \end{enumerate}
\end{prop}
\begin{pf}
(1)$\Rightarrow$(2)$\Rightarrow$(3) is obvious. To show that 
(3)$\Rightarrow$(1) let $(x_n)$ be a dense sequence $(x_n)$ in the 
dense subset of $A$ given in (3) and  let $(t_n)$ be a sequence in 
$\R$; then there is a subsequence $(t_{n(k)})$ such that both 
$(\alpha_{t_{n(k)}}(x_m))$ and  $(\alpha_{-t_{n(k)}}(x_m))$ 
converge for any $m$. Then the (strong) limit of 
$(\alpha_{t_{n(k)}})$ exists as an automorphism. This implies (1). 

If (1) holds, then the closure of $\alpha_{\bf R}$ is a compact 
abelian group and hence (4) follows. It is immediate that 
(4)$\Rightarrow$(5)$\Rightarrow$(3). 
\end{pf} 

We recall that $\alpha$ is said to be a compact flow if the 
conditions in \ref{D1} are satisfied. 
                           
\ncm{\PSp}{{\rm PSp}}
\begin{prop}\label{D2}
Let $A$ be a unital separable simple \cstar\ and $\alpha$ a flow 
on $A$. If $\alpha$ is an approximately inner compact flow, then 
$D(\delta_{\alpha})$ contains a maximal abelian $C^*$-subalgebra 
of $A$. \end{prop} 
\begin{pf}
Since $\alpha$ is approximately inner, $\alpha$ has a pure ground 
state $\omega$. If $G$ denotes the closure of $\alpha_{\bf R}$, 
then $\omega$ is left invariant under $G$. Thus the GNS 
representation associated with $\omega$ is a $G$-covariant 
irreducible representation $\pi$ of $A$. Since $A$ is simple, 
$\pi$ is faithful.     

Let $U$ be a continuous representation of $G$ on ${\cal H}_{\pi}$ 
such that $U_g\pi(x)U_g^*=\pi g(x)$ for $x\in A$. There is an 
orthogonal family $(E_p)_{p\in \hat{G}}$ of projections such that 
$U_g=\sum_p<g,p>E_p$, where $\hat{G}$ is the character group of 
$G$; $\hat{G}$ is a countable discrete abelian group.  Let 
$A^G=\{x\in A\ |\ \forall g\in G,\, g(x)=x\}$ and denote by 
$\pi_p$ the representation of $A^G$ on $E_p{\cal H}_{\pi}$, i.e., 
$\pi_p(a)=\pi(a)|E_p{\cal H}_{\pi},\ a\in A^G$. It follows that 
$\{\pi_p\,|\,p\in\hat{G}\}$ is a disjoint family of irreducible 
representations of $A^G$ since $\pi(A^G)'=U_G''$. Under this 
circumstance we will show that there is a maximal abelian 
$C^*$-subalgebra (masa) of $A^G$ which is also a masa in $A$. 
   
By the following lemma there is an $h\in (A^G)_{sa}$ such that 
$\pi(h)\cong\oplus_{p\in\hat{G}} \pi_p(h)$ is diagonal and all the 
eigenvalues of $\pi(h)$ have multiplicity one.  Let $C$ be a masa 
of $A^G$ such that $C\ni h$. Since $A\cap C'$ is $G$-invariant, 
$A\cap C'$ is the closed linear span of the eigenspaces of 
$G|A\cap C'$. Let $x\in A\cap C'$ be a nonzero element such that 
$g(x)=<g,q>x,\ g\in G$ for some $q\in \hat{G}$. Then there must be 
a $p\in\hat{G}$ and a unit vector $\xi$ in $E_p{\cal H}$  such 
that $\pi(x)\xi\neq0$ and $\pi(h)\xi=\lambda\xi$ for some 
$\lambda\in \R$. Since $\pi(x)\xi\in E_{p+q}{\cal H}$ and 
$\pi(h)\pi(x)\xi=\pi(x)\pi(h)\xi=\lambda\pi(x)\xi$, $\pi(x)\xi$ is 
a constant multiple of $\xi$. This implies that $q=0$ and then 
$x\in C$. Thus we can conclude that $C$ is a masa in $A$.  
\end{pf} 
                                               
\begin{lem}\label{D3}
Let $A$ be a separable $C^*$-algebra and let $(\pi_n)$ be a 
sequence of irreducible representations of $A$ such that if $m\neq 
n$ then $\pi_m$ is disjoint from $\pi_n$. If $(S_n)$ is a sequence 
of dense subsets of $\R$ then there is an $h\in A_{sa}$ such that 
$\pi_n(h)$ is diagonal and $\PSp(\pi_n(h))\subset S_n$ for all 
$n$.  Furthermore $h$ can be chosen so that all the eigenvalues of 
$\pi_n(h)$ have multiplicity one for every $n$. 
\end{lem} 
\begin{pf}
We will construct a sequence $(h_n)$ in $A_{sa}$ and an increasing 
sequence $(e_{kn})_{n\geq k}$ of finite-dimensional projections on 
${\cal H}_k={\cal H}_{\pi_k}$ for each $k$ such that 
$\|h_n\|<2^{-n}$, $\lim_n e_{kn}=1$, $\pi_k(h_n)e_{k,n-1}=0$ for 
$k=1,2,\ldots,n-1$, and for $H_n=h_1+h_2+\cdots+h_n$, 
 \BE                             
 &&[\pi_k(H_n),e_{kn}]=0,\ \ k=1,2,\ldots,n\\
 && \PSp(\pi_k(H_n)|e_{kn}{\cal H}_k)\subset S_k, \ \ 
 k=1,2,\ldots,n,
 \EE            
where $\PSp$ denotes the set of eigenvalues. Then we will let 
$h=\lim_nH_n$ and then since $\pi_k(h)e_{k,n}=\pi_k(H_n)e_{kn}$, 
we will have that $\PSp(\pi_k(h))\subset S_k$ for all $k$. 
        
To construct such sequences as above, we will argue inductively by 
using Kadison's transitivity theorem (see, e.g., 1.21.16 of 
\cite{Sakai}). Suppose that we have constructed 
$h_1,h_2,\ldots,h_n$ and $e_{km}$ with $n\geq m\geq k$ satisfying 
the above conditions. With $H_n=h_1+\cdots h_n$ let $E_k$ be the 
spectral measure of $\pi_k(H_n)$ for $k=1,2,\ldots,n+1$. For each 
$k$ we find a finite family $F_k$ of disjoint translates of 
$[0,2^{-n-1})$ and a family $\{\xi(I)\}_{I\in F_k}$ of unit 
vectors in $(1-e_{kn}){\cal H}_k$ such that $E_k(I)\xi(I)=\xi(I)$ 
and the linear span of $\xi(I),\,I\in F_k$ is so large that it 
almost contains any prescribed vector from $(1-e_{kn}){\cal H}_k$. 
Choose $\lambda_I\in S_k\cap I$ and let $P_I$ denote the  
projection onto the linear space spanned by $\xi(I)$ and 
$(\pi_k(H_n)-\lambda_I)\xi(I)$. Then the family $(P_I)_{I\in F_k}$ 
is mutually orthogonal and the sum $P_k=\sum_IP_I$ is orthogonal 
to $e_{kn}$. It follows that 
$\|P_k(\pi_k(H_n)-\sum_I\lambda_IP_I)P_k\|= 
\|\sum_IP_I(\pi_k(H_n)-\lambda_IP_I)P_I\|<2^{-n-1}$. We will then 
choose an $h_{n+1}\in A_{sa}$ such that $\|h_{n+1}\|\leq2^{-n-1}$,  
$[\pi_k(h_{n+1}),P_k]=0$ for $k\leq n+1$, $\pi_k(h_{n+1})e_{kn}=0$ 
for $k\leq n$, and 
$\pi_k(h_{n+1})P_k=-(\pi_k(H_n)-\sum_I\lambda_IP_I)P_k$ for 
$k=1,2,\ldots,n+1$. Then we obtain that 
$\pi_k(H_n+h_{n+1})\xi(I)=\lambda(I)\xi(I)$. We let $e_{k,n+1}$ be 
the sum of $e_{kn}$ (0 if $k=n+1$) and the projection onto the 
linear span of $\xi(I),\ I\in F_k$. Thus we have constructed 
$h_{n+1},\, e_{k,n+1}$ as required. (A similar argument is used in 
\cite{K88}.) 

\end{pf}  
  
\begin{rem}\label{D4}
In \ref{D2} the assumption that $\alpha$ is approximately inner is 
made to ensure that there is a $G$-covariant irreducible 
representation for $G=\ol{\alpha_{\bf R}}$. This assumption 
certainly is not necessary as we can see from the following 
example.  If $\alpha$ is a flow on the Cuntz algebra ${\cal 
O}_2=C^*(s_1,s_2)$ such that $\alpha_t(s_j)=e^{i\mu_j t}s_j$, 
where $\mu_1,\mu_2$ are rationally independent, then $\alpha$ is 
compact and $G$ is the gauge action $\gamma$ of $\T^2$ given by 
$\gamma_{z_1,z_2}(s_j)=z_js_j$. In this case $\alpha$ is not 
approximately inner but $\gamma$ has a covariant irreducible 
representation. Hence the conclusion of \ref{D2} follows for this 
$\alpha$ (though it is easy to show that directly). On the other 
hand there is an example of a compact flow where $G$ does not have 
a covariant irreducible representation. For example if $\alpha$ is 
a flow on the irrational rotation \cstar\ 
$A_{\theta}=C^*(u_1,u_2)$ with $u_1u_2=e^{2\pi i\theta}u_2u_1$ and 
$\theta$ irrational such that $\alpha_t(u_j)=e^{i\mu_j t}u_j$, 
where $\mu_1,\mu_2$ are rationally independent, then $G$ is the 
gauge action $\gamma$ of $\T^2$ given by 
$\gamma_{z_1,z_2}(u_j)=z_ju_j$. Since $\gamma$ is ergodic, only 
the tracial representation is covarint though there is an 
$\alpha$-covariant irreducible representation \cite{K88}. In this 
case we do not know whether the domain $D(\delta_{\alpha})$ 
contains a masa or not. 
\end{rem}

\section{Universal UHF flows}
 \setcounter{theo}{0} 
We say that a flow $\alpha$ on $A$ is {\em almost conjugate} to a 
flow $\beta$ on $B$ and denote it by 
$\alpha\stackrel{ac}{\sim}\beta$ if for any $\epsilon>0$ there is 
an isomorphism $\varphi$ of $A$ onto $B$ such that 
$\|\alpha_t-\varphi^{-1}\beta_t\varphi\|<\epsilon$ for $t\in 
[-1,1]$. If $A$ (or $B$) is separable and simple and 
$\alpha\stackrel{ac}{\sim}\beta$, then 
$\alpha\stackrel{cc}{\sim}\beta$, i.e., $\alpha$ is cocycle 
conjugate to $\beta$ (see 1.2 of \cite{K5}). To prove this we use 
the fact that if $A$ is simple and 
$\|\alpha_t-\varphi^{-1}\beta_t\varphi\|<2$, then there is a 
unitary $u_t\in A$ such that 
$\alpha_t=\Ad\,u_t\varphi^{-1}\beta_t\varphi$. To get an 
$\alpha$-cocycle from the family $(u_t)$ of unitaries so obtained 
for small $|t|$, we may use 8.1 of \cite{OL}, which requires 
separability.   
 
Let $A$ be a UHF \cstar\ such that $A\otimes A\cong A$. We recall 
that a flow $\gamma$ on $A$ is a universal UHF flow if 
$\gamma\otimes\alpha \stackrel{cc}{\sim}\gamma$ for any UHF flow 
$\alpha$ on $A$. 

\begin{prop}\label{A1} 
If $A$ is a UHF \cstar\ such that $A\otimes A\cong A$, then there 
is a universal UHF flow $\gamma$ on $A$. 
\end{prop} 
\begin{pf}
There exists a (finite or infinite) sequence $(p_i)$ of prime 
numbers such that $A\cong \otimes_iM_{p_i^{\infty}}$. Let $S$ be 
the set of integers which are of the form 
$p_1^{k_i}p_2^{k_2}\cdots p_n^{k_n}$ with $k_i\geq0$. Then it 
follows that $A\cong \otimes_{q\in S}M_{q^{\infty}}$. For each 
$q\in S$ let $(h_n)$ be a dense sequence in the self-adjoint 
diagonal matrices of $M_q$.  We define a flow $\gamma^{(q)}$ on 
the UHF \cstar\ $M_{q^{\infty}}$ by 
 $$
 \bigotimes_{n=1}^{\infty}\Ad\,e^{ith_n}
 $$
and define a flow $\gamma$ on $A$ by $\otimes_{q\in 
S}\gamma^{(q)}$. 

If $\alpha$ is a UHF flow on $A$, then there are an infinite 
sequence $(q_n)$ in $S$ and a sequence $(k_n)$ with $k_n\in 
M_{q_n}$ being self-adjoint and diagonal such that $\alpha$ is 
conjugate to 
 $$
 \bigotimes_{n=1}^{\infty}\Ad\,e^{itk_n}.
 $$
Hence $\gamma\otimes\alpha$ is of the same form as 
$\gamma=\otimes_{q\in S}\gamma^{(q)}$ (by tensoring $\gamma^{(q)}$ 
with $\otimes_{n:q_n=q}\Ad\,e^{itk_n}$) but the sequence $(h_n')$ 
defining the flow on $M_{q^{\infty}}$ may be different from the 
$(h_n)$ defining $\gamma^{(q)}$ for some $q\in S$. But since they 
are dense in the self-adjoint diagonal matrices of 
$M_{q^{\infty}}$ for any $q\in S$, we can conclude that 
$\gamma\otimes\alpha$ is cocycle conjugate to $\gamma$, i.e., 
$\gamma\otimes\alpha\stackrel{cc}{\sim}\gamma$. (As a matter of 
fact $\gamma\otimes\alpha\stackrel{ac}{\sim}\gamma$.)  
\end{pf}
  
Let $F(A)$ denote the set of flows on $A$. We give a topology on 
$F(A)$ by
 $$
 d(\alpha,\beta)=\sum_{n=1}^{\infty}2^{-n}
 \max_{|t|\leq1}\|\alpha_t(x_n)-\beta_t(x_n)\|,
 $$
where $(x_n)$ is a dense sequence of the unit ball of $A$. That 
is, $\alpha_n$ converges to $\alpha$ if $\alpha_{nt}(x)$ converges 
to $\alpha_t(x)$ uniformly in $t$ on $[-1,1]$ for all $x\in A$. 
Let $AIF(A)$ denote the set of approximately inner flows on $A$, 
i.e., the closure of inner flows.

\begin{prop}\label{A2}
Let $A$ be a UHF \cstar\ such that $A\otimes A\cong A$ and 
$\gamma$ a universal UHF flow on $A$. Then
 $$\{\varphi^{-1}\gamma\varphi\ |\ \varphi\in \Aut(A)\}\equiv C(\gamma)
 $$
is dense in $AIF(A)$.  
\end{prop}
\begin{pf}
Let $h\in A_{sa}$. It suffices to show that the closure of 
$C(\gamma)$ contains the inner flow $t\mapsto\Ad\,e^{ith}$. 

Since $\gamma$ is a universal UHF flow, it follows that 
$\gamma\otimes(\otimes^{\infty}\Ad\,e^{ith}) 
\stackrel{cc}{\sim}\gamma$, i.e., there is an isomorphism 
$\varphi$ of $A\otimes(\otimes^{\infty}A)\equiv \otimes^{\infty}A$
onto $A$ and a $k\in A_{sa}$ such that 
 $$
 \gamma_t\otimes(\bigotimes^{\infty}\Ad\,e^{ith})
 =\varphi^{-1}\gamma_t^{(k)}\varphi,
 $$                                 
where $\gamma^{(k)}$ is the inner perturbation of $\gamma$ by 
$\ad\,ik$.  Let $\varphi_n$ denote the homomorphism of $A$ into 
$A$ obtained by restricting $\varphi$ to the $n+1$'st factor of 
$\otimes^{\infty}A$. Then we have that 
 $$
 \lim_{k\rightarrow\infty}(\varphi_n\Ad\,e^{ith}(x)-\gamma_t\varphi_n(x))=0
 $$
uniformly in $t$ on $[-1,1]$ for every $x\in A$. Since there are 
isomorphisms $\tilde{\varphi}_n$ of $A$ onto $A$ such that 
$\lim(\varphi_n(x)-\tilde{\varphi}_n(x))=0$ for every $x\in A$, we 
may assume that $\varphi_n$'s are all isomorphisms of $A$ onto $A$ 
and obtain that $\Ad\,e^{ith}$ is the limit of 
$\varphi_n^{-1}\gamma\varphi_n\in C(\gamma)$.
\end{pf}

\begin{rem}\label{A3}
If a flow $\alpha$ absorbs the universal UHF flow $\gamma$, i.e., 
$\alpha\otimes\gamma\stackrel{cc}{\sim}\alpha$, then the closure 
of $C(\alpha)$ contains $AIF(A)$. This follows since for $h\in 
A_{sa}$,
 $$
 \alpha\otimes(\bigotimes^{\infty}\Ad\,e^{ith})\stackrel{cc}{\sim}
 \alpha\otimes\gamma\otimes(\bigotimes^{\infty}\Ad\,e^{ith})\stackrel{cc}{\sim} 
 \alpha\otimes\gamma \stackrel{cc}{\sim}\alpha.
 $$
\end{rem}

\begin{prop}\label{A4}
If $A$ is a UHF \cstar\ such that $A\otimes A\cong A$, then the 
universal UHF flows on $A$ are mutually almost conjugate.  In 
particular a cocycle perturbation of a universal UHF flow $\gamma$ 
is almost conjugate to $\gamma$. 
\end{prop}
\begin{pf}
If $\alpha$ is a UHF flow and $\gamma$ is a universal UHF flow, it 
follows that 
$\alpha^{\infty}\otimes\gamma^{\infty}\stackrel{cc}{\sim}\gamma$, 
where $\alpha^{\infty}=\otimes^{\infty}\alpha$ etc. Hence for any 
$\epsilon>0$ there exists a flow $\rho$ on $A$ such that 
$\alpha^{\infty}\otimes\gamma^{\infty}\otimes\rho$ is conjugate to 
$\gamma$ up to $\epsilon$ (i.e., 
$\|\alpha^{\infty}_t\otimes\gamma^{\infty}_t\otimes\rho_t 
-\varphi^{-1}\gamma_t\varphi\|<\epsilon$ for any $t\in [-1,1]$ and 
for some isomorphism $\varphi:A^{\infty}\ra A$). Hence 
$\alpha^{\infty}\otimes\gamma^{\infty}\otimes\rho=\alpha\otimes 
\alpha^{\infty}\otimes\gamma^{\infty}\otimes\rho$ is conjugate to 
$\alpha\otimes\gamma$ up to $\epsilon$. Hence 
$\alpha\otimes\gamma\stackrel{ac}{\sim}\gamma$. If both $\alpha$ 
and $\gamma$ are universal, it follows that 
$\alpha\stackrel{ac}{\sim}\alpha\otimes\gamma\stackrel{ac}{\sim} 
\gamma$. 
\end{pf} 

\begin{rem}\label{A5}
Let $A$ be a UHF \cstar\ with $A\otimes A\cong A$. Then there is a 
{\em universal} automorphism $\gamma$ of $A$ in the sense that 
$\gamma\otimes\alpha$ is cocycle conjugate to $\gamma$ for any 
automorphism $\alpha$ of $A$. Namely an automorphism with the 
Rohlin property is universal (see, e.g., \cite{BKRS,K1,K0,EK}). 
The universal automorphisms are mutually almost conjugate. 
\end{rem} 

\section{The flip automorphism}
 \setcounter{theo}{0} 
When $A$ is a UHF \cstar, the flip automorphism $\sigma$ of 
$A\otimes A$ defined by $\sigma(x\otimes y)=y\otimes x,\ x,y\in A$ 
is approximately inner (since $\sigma$ induces the identity map on 
$K_0(A\otimes A)$, which has rank one). Hence there is a sequence 
$(u_n)$ of unitaries in $A\otimes A$ such that 
$\sigma=\lim_n\Ad\,u_n$. 

Let $\alpha$ be a flow on $A$. Since $\sigma(\alpha\otimes\alpha)= 
(\alpha\otimes\alpha)\sigma$, we can ask a question of whether one 
can put an extra condition, in connection with $\alpha$, on the 
sequence $(u_n)$ above. 

Before turning to this special situation we present the following 
two propositions to clarify the conditions we are thinking of (cf. 
\cite{FVV}).       

\begin{prop}\label{B1}  
Let $A$ be a unital \cstar\ and let $\alpha$ be a flow on $A$ and 
$\sigma\in \Aut(A)$ such that $\sigma^{-1}\alpha\sigma=\alpha$. 
Then the following conditions are equivalent: 
\begin{enumerate}
\item There is a sequence $(u_n)$ of unitaries in 
$D(\delta_{\alpha})$ such that $\sigma=\lim_n\Ad\,u_n$ and 
$\sup_n\|\delta_{\alpha}(u_n)\|<\infty$. 
 \item There is a sequence $(u_n)$ of unitaries in $A$ 
 such that $\sigma=\lim_n\Ad\,u_n$ and $(t\mapsto 
 \alpha_t(u_n))$ are equi-continuous in $n$.
 \end{enumerate}
\end{prop}
\begin{pf}
If (1) is satisfied, then so is (2) for the same $(u_n)$. To go 
from (2) to (1) we may need to modify $(u_n)$. Since a similar 
technique will apply in the proof of the following proposition, we 
do not give it here. \end{pf}

We will express the conditions in the above proposition by saying 
that $\sigma$ is {\em $\alpha$-smoothly approximately inner}, 
which is weaker than the condition that $\sigma$ is {\em 
$\alpha$-invariantly approximately inner}, which will appear in 
the following.  

\begin{prop}\label{B2} 
Let $A$ be a unital \cstar\ and let $\alpha$ be a flow on $A$ and 
$\sigma\in\Aut(A)$ such that $\sigma^{-1}\alpha\sigma=\alpha$. 
Then the following conditions are equivalent: 
\begin{enumerate}
\item There is a sequence $(u_n)$ of unitaries in 
 $D(\delta_{\alpha})$ such that 
 $\sigma=\lim_n\Ad\,u_n$ and 
 $\lim_n\|\delta_{\alpha}(u_n)\|=0$. 
 \item There is a sequence $(u_n)$ of unitaries in $A$ 
 such that $\sigma=\lim_n\Ad\,u_n$ and 
 $\|\alpha_t(u_n) -u_n\|\ra0 $ uniformly in $t$ on
 every compact subset of $\R$.
 \item There is a sequence $(u_n)$ of unitaries in $A$ 
 such that $\sigma=\lim_n\Ad\,u_n$ and 
 $\|\alpha_t(u_n) -u_n\|\ra0 $ for every $t\in\R$.
 \end{enumerate}
\end{prop} 
\begin{pf}
It is obvious that (1)$\Rightarrow$(2)$\Rightarrow$(3).

(3)$\Rightarrow$(2).  Let $\epsilon>0$ and define
 $$
 I_n=\{t\in \R\ |\ 
 \sup_{m\geq n}\|\alpha_t(u_m)-u_m\|\leq\epsilon\}.  
 $$
Since $\cup_nI_n=\R$ and $I_n$'s are closed, $\cup_nI_n^{\circ}$ 
is dense in $\R$ by the Baire Category theorem, where 
$I_n^{\circ}$ is the interior of $I_n$. If $(a,a+\delta)\subset 
I_n^{\circ}$ for some $a\in \R$, $\delta>0$, and $n$, then 
$\|\alpha_t(u_m)-u_m\|\leq 2\epsilon$ for $t\in (-\delta,\delta)$ 
and $m\geq n$. Then the rest is easy. 

(2)$\Rightarrow$(1). Let $f$ be a $C^{\infty}$-function on $\R$ 
with compact support such that $f\geq0$ and  $\int f(t)dt=1$ and 
let, for a small $\epsilon>0$, 
 $$
 x_{n\epsilon}=\epsilon\int f(\epsilon t)\alpha_t(u_n)dt.
 $$
Since 
 $$
 \|x_{n\epsilon}-u_n\|\leq 
 \epsilon\int f(\epsilon t)\|\alpha_t(u_n)-u_n\|dt,
 $$
one can choose a sequence $(\epsilon_n)$ such that 
$\epsilon_n\searrow0$ and $\|y_n-u_n\|\ra0$ with 
$y_n=x_{n\epsilon_n}$. Note that $y_n\in D(\delta_{\alpha})$ and 
that 
 $$
 \|\delta_{\alpha}(y_n)\|\leq \epsilon_n\int|f'(t)|dt.
 $$ 
If we denote by $v_n$ the unitary obtained by the polar 
decomposition of $y_n$, the sequence $(v_n)$ satisfies the desired 
properties. \end{pf}

If a flow $\alpha$ has more than one $\alpha$-KMS states for some 
inverse temperature, the following shows that the flip is not 
$\alpha$-smoothly approximately inner. 
 
\begin{prop}\label{B3}
Let $A$ be a UHF \cstar\ and $\alpha$ a flow on $A$ and suppose 
that the flip is $\alpha$-smoothly approximately inner (or more 
precisely, $\alpha\otimes\alpha$-smoothly approximately inner), 
i.e., the conditions in \ref{B1} are satisfied for $A\otimes A$, 
$\alpha\otimes\alpha$, and the flip automorphism $\sigma$ in place 
of $A$, $\alpha$, and $\sigma$ respectively. Then the the set of 
$\alpha$-KMS states is a singleton (if not empty) for each inverse
temperature. 
\end{prop} 
\begin{pf}
Let $\omega_i$ be an $\alpha$-KMS state at inverse temperature 
$c\in\R$ for $i=1,2$. Then $\omega_1\otimes\omega_2$ is an 
$\alpha\otimes\alpha$-KMS state of $A\otimes A$ at  $c$. Since the 
flip $\sigma$ commutes with $\alpha\otimes\alpha$ and the flip is 
$\alpha$-smoothly approximately inner, the flip leaves each 
$\alpha\otimes\alpha$-KMS state invariant by the following lemma, 
i.e., $(\omega_1\otimes\omega_2)\sigma=\omega_1\otimes\omega_2$ or 
$\omega_2\otimes\omega_1=\omega_1\otimes\omega_2$. Hence 
$\omega_1=\omega_2$. \end{pf}                                   

The following lemma is shown by M. Fannes et al \cite{FVV} (or 
5.3.33A of \cite{BR} vol.~II), but we will present another proof 
based on the definition of the KMS condition in terms of 
holomorphic functions. 

\begin{lem} \label{B35}
Let $A$ be a unital \cstar\ and let $\alpha$ be a flow and 
$\sigma\in\Aut(A)$ such that $\sigma^{-1}\alpha\sigma=\alpha$. 
Suppose that $\sigma$ is $\alpha$-smoothly approximately inner. 
Then $\omega\sigma=\omega$ for any $\alpha$-KMS states at inverse 
temperature $c\in\R$. 
\end{lem} 
\begin{pf}
Let $\omega$ be an $\alpha$-KMS state at inverse temperature 
$c>0$. (The case $c<0$ entails just some notational changes and 
the case $c=0$ is trivial.) By definition, for any $x,y\in A$ 
there is a bounded continuous function $F(z)$ on $S_c=\{z\in \C\ 
|\ 0\leq\Im z\leq c\}$ such that $F$ is analytic in the interior 
of $S_c$ and 
 $$
 F(t)=\omega(x\alpha_t(y)),\ \ F(t+ic)=\omega(\alpha_t(y)x)
 $$                                                        
for $t\in\R$. Let $(u_n)$ be the sequence of unitaries as given in 
\ref{B1} and denote by $F_{n,x}$ the function $F$ obtained by 
taking $u_nx$ for $x$ and $u_n^*$ for $y$. In particular $F_{n,x}$ 
satisfies 
 $$
 F_{n,x}(t)=\omega(u_nx\alpha_t(u_n^*)),\ \ 
 F_{n,x}(t+ic)=\omega(\alpha_t(u_n^*)u_nx),
 $$                                        
for $t\in\R$. Since $t\mapsto \alpha_t(u_n^*)u_n$  is 
equi-continuous in $n$, we may assume by passing to a subsequence 
that $F_{n,y}$ converges, for $y=x$ and $y=1$, uniformly on every 
compact subset of the boundary $\partial S_c$. Then $F_{n,y}$ 
converges, say to $F_y$, uniformly on every compact subset of 
$S_c$. It follows that $F_y$ is a bounded continuous function on 
$S_c$ which is analytic in the interior. 

  Suppose further that $\omega$ is factorial, which causes no 
loss of generality since the extreme KMS states are all factorial. 
Then since $(\alpha_t(u_n^*)u_n) $ is central, it follows that 
$\omega(\alpha_t(u_n^*)u_nx)$ converges to $F_1(t+ic)\omega(x)$. 
This implies that $F_x(z)=F_1(z)\omega(x)$ for $z\in\R+ic$ and 
hence for $z\in S_c$. Since $F_x(0)=\omega\sigma(x)$ and 
$F_1(0)=1$, we obtain that $\omega\sigma(x)=\omega(x)$. Hence 
$\omega\sigma=\omega$. 

\end{pf}   

We may apply the above proposition \ref{B3} to the one-dimensional 
quantum lattice systems (since the {\em bounded surface energy} 
condition obviously implies that the flip is $\alpha$-smoothly 
approximately inner), where the uniqueness of KMS states is of 
course well-known (see \cite{Sak}). This is yet another proof.
      
When $L$ is a bounded linear map of a \cstar\ $B$ into a \cstar\ 
$A$, we denote by $\|L\|_{cb}$ the completely bounded norm defined 
by $\sup_n\|L\otimes \id_n\|$, where $\id_n$ is the identity map 
on the matrix algebra $M_n$ so that $L\otimes\id_n$ is a linear 
map of $B\otimes M_n$ into $A\otimes M_n$. In the following 
proposition we will use such a norm when $B$ is finite-dimensional 
and $L$ is a difference of homomorphisms. In this case it follows 
from Christensen's result \cite{Ch} that for any $\epsilon>0$ 
there is a $\delta>0$ such that if $\|L\|<\delta$ then 
$\|L\|_{cb}<\epsilon$. (Because if $L=\phi_1-\phi_2$ has 
sufficiently small norm  with $\phi_i$ a unital homomorphism of 
$B$ into $A$, then there is a unitary $u\in A$ such that $\|u-1\|$ 
is small and $\phi_2=\Ad\,u\phi_1$, which entails that 
$\|L\|_{cb}=\|(\id-\Ad\,u)\phi_1\|_{cb}\leq 2\|u-1\|$.)

\begin{prop}\label{B4}
Let $A$ be a UHF \cstar\ with $A\otimes A\cong A$ and $\alpha$ a 
flow on $A$. Suppose that there exists a sequence $(\varphi_n)$ in 
$\Aut(A)$ such that $\varphi_n^{-1}\alpha\varphi_n\ra \alpha$ and 
$(\varphi_n(x))$ is central for any $x\in A$ and that   the flip 
is $\alpha$-invariantly approximately inner. 

Then for any unital finite type I subfactor $A_1$ of $A$ and 
$\epsilon>0$ there is a finite type I subfactor $B$ of $A$ with 
$B\supset A_1$ such that if $\varphi\in \Aut(A)$ satisfies that 
for any  $t\in I=[-1,1]$, 
 $$
 \|(\alpha_t-\varphi^{-1}\alpha_t\varphi)|B\|_{cb}< \epsilon'
 $$
for some $\epsilon'>0$, then  there is a unitary $u\in A$ such 
that 
 \BE
 \max_{t\in I}\|\alpha_t(u)-u\|&<&\epsilon'+\epsilon,\\ 
 \|(\varphi-\Ad\,u)|A_1\|&<& \epsilon.
 \EE
\end{prop}
\begin{pf}
Let $\sigma$ denote the flip of $A\otimes A$. By the assumption on 
$\alpha$ there exists a unitary $u\in A\otimes A$ such that
 \BE
 \Ad\,u|A_1\otimes A_1 &=& \sigma|A_1\otimes A_1,\\
 \max_{t\in I}\|\alpha_t\otimes\alpha_t(u)-u\|&<&\epsilon/2.
 \EE
We may suppose that there is a finite type I subfactor $B$ of $A$ 
such that $B\otimes B\ni u$ and $B\supset A_1$.  

By using the $(\varphi_n)$ and $\varphi$ in the statement, we 
define linear maps $\iota\otimes\varphi_n$ and 
$\varphi\otimes\varphi_n$ of the algebraic tensor product $A\odot 
A$ into $A$ by
 \BE
  \iota\otimes\varphi_n(x\otimes y)&=& x\varphi_n(y),\\
  \varphi\otimes\varphi_n(x\otimes y)&=&\varphi(x)\varphi_n(y).
  \EE
Since $(\varphi_n(y))$ is a central sequence for any $y\in A$, 
both  $\iota\otimes\varphi_n$ and $\varphi\otimes\varphi_n$ are 
{\em approximate} homomorphisms. Thus $\iota\otimes\varphi_n(u)$ 
and $\varphi\otimes\varphi_n(u)$ are close to unitaries for the 
$u\in B\otimes B$ above for all large $n$.

Since 
 $$\|\alpha_t(\iota\otimes\varphi_n)(u)
 -\iota\otimes\varphi_n(\alpha_t\otimes\alpha_t(u))\|
 $$
converges to zero uniformly in $t\in I$, we obtain that  
 $$
 \| \alpha_t(\iota\otimes\varphi_n)(u)
 -\iota\otimes\varphi_n(u)\|<\epsilon/2
 $$
for $t\in I$ and for all large $n$.

For a finite sum $\sum_ix_i\otimes y_i\in B\otimes B$ it follows 
that 
 $$
 \|\alpha_t( \varphi\otimes\varphi_n)(\sum_ix_i\otimes y_i)
 -( 
 \varphi\otimes\varphi_n)(\sum_i\alpha_t(x_i)\otimes\alpha_t(y_i))\|
 $$
converges to 
 $$
 \|\sum_i\alpha_t\varphi(x_i)\otimes\alpha_t(y_i)
 -\sum_i\varphi\alpha_t(x_i)\otimes\alpha_t(y_i)\|
 $$
which is less than or equal to $ 
 \|(\alpha_t\varphi-\varphi\alpha_t)|B\|_{cb}\,\|\sum_ix_i\otimes y_i\|$.
Hence it follows that
 $$
 \|\alpha_t( \varphi\otimes\varphi_n)(u) -( 
 \varphi\otimes\varphi_n)(\alpha_t\otimes\alpha_t(u))\|<\epsilon'
 $$
for any $t\in I$ and for all large $n$ and that
 $$
 \max_{t\in I}\|\alpha_t( \varphi\otimes\varphi_n)(u) -( 
 \varphi\otimes\varphi_n)(u)\|<\epsilon'+\epsilon/2.
 $$
Thus we obtain that $\max_{t\in 
I}\|\alpha_t(y_n)-y_n\|<\epsilon'+\epsilon$ for 
$y_n=(\varphi\otimes\varphi_n)(u)(\iota\otimes\varphi_n)(u)$.

Since, for $x\in A_1$,
 $$
 (\iota\otimes\varphi_n)(u)x\approx ( 
 \iota\otimes\varphi_n)(u\cdot x\otimes 1)= 
 (\iota\otimes\varphi_n)(1\otimes x\cdot u)\approx 
 \varphi_n(x)(\iota\otimes\varphi_n)(u)
 $$
and
 $$
 (\varphi\otimes\varphi_n)(u)\varphi_n(x)\approx
   (\varphi\otimes\varphi_n)(u\cdot 1\otimes x)=
   (\varphi\otimes\varphi_n)(x\otimes 1\cdot u)\approx
   \varphi(x)(\varphi\otimes\varphi_n)(u),
 $$
it follows that $\|y_nx-\varphi(x)y_n\|\ra0$ for any $x\in A_1$ as 
$n\ra\infty$. Thus, by choosing a sufficiently large $n$ and 
taking the unitary part of the polar decomposition of $y_n$, which 
is already close to a unitary, we obtain the conclusion. (As a 
matter of fact $\|(\Ad\,u-\varphi)|A_1\|$ can be made small 
independently of $\epsilon$.) 
\end{pf}  
 
\begin{theo}\label{B5}
Let $A$ be a UHF \cstar\ with $A\otimes A\cong A$ and $\gamma$ a 
universal UHF flow on $A$. If $\alpha$ is an approximately inner 
flow on $A$ such that 
$\alpha\otimes\gamma\stackrel{cc}{\sim}\alpha$, the following 
conditions are equivalent: 
 \begin{enumerate}
 \item $\alpha\stackrel{cc}{\sim}\gamma$;     
 \item The flip is $\alpha$-invariantly approximately inner.
 \end{enumerate}
\end{theo}
\begin{pf}
(1)$\Rightarrow$(2). We may assume that $\alpha$ is an inner 
perturbation of $\gamma$, i.e., 
$\delta_{\alpha}=\delta_{\gamma}+\ad\,ih$ for some $h\in A_{sa}$. 
We know that there is a sequence $(u_n)$ of unitaries in 
$D(\delta_{\gamma})\odot D(\delta_{\gamma})$ such that 
$u_n^*=u_n$, $\delta_{\gamma\otimes\gamma}(u_n)=0$, and 
$\sigma=\lim_n\Ad\,u_n$, where $\sigma$ is the flip of $A\otimes 
A$. Hence it follows that $u_n\in D(\delta_{\alpha\otimes\alpha})$ 
and that
 $$
 \delta_{\alpha\otimes\alpha}(u_n)=\ad(ih\otimes1+1\otimes ih)(u_n)
 =(ih\otimes1)u_n-u_n(1\otimes ih)+(1\otimes ih)u_n-u_n(ih\otimes1)
 $$
converges to zero. This is what we wanted to show; see \ref{B2}. 

(2)$\Rightarrow$(1).  Since 
$\gamma^{\infty}\stackrel{cc}{\sim}\gamma$, we see that there is a 
sequence $(\varphi_n)$ in $\Aut(A)$ such that 
$\varphi_n^{-1}\gamma\varphi_n\ra \gamma$ and $(\varphi_n(x))$ is 
central for any $x\in A$. Since $\alpha$ is an approximately inner 
flow satisfying 
$\alpha\stackrel{cc}{\sim}\alpha\otimes\gamma^{\infty}$ and any 
inner flow can be embedded into $\gamma$, we see that $\alpha$ 
also satisfies the above condition. Since the flip is both 
$\gamma$ and $\alpha$-invariantly approximately inner, we can 
apply \ref{B4} to both of them.

Let $(A_n)$ be an increasing sequence of finite type I subfactors 
of $A$ with $A=\ol{\cup_nA_n}$.

Let $\epsilon>0$. For $A_1$, $\alpha$, and $2^{-1}\epsilon$ (in 
place of $\epsilon$) we choose $B_1$ as $B$ in \ref{B4}. By 
\ref{A2} we choose $\varphi_1\in \Aut(A)$ such that for $t\in 
I=[-1,1]$, 
 $$
 \|(\alpha_t-\varphi_1^{-1}\gamma_t\varphi_1)|B_1\|_{cb}<2^{-2}\epsilon.
 $$
By slightly changing $\varphi_1$ if necessary we assume that 
$\varphi_1(B_1)\subset A_n$ for some $n>1$; by passing to a 
subsequence of $(A_n)$ we assume that $n=2$; i.e., 
$\varphi_1(B_1)\subset A_2$. We then choose $B_2$ for $A_2$, 
$\gamma$, and $2^{-2}\epsilon$ by \ref{B4} and choose 
$\varphi_2\in \Aut(A)$ by \ref{A3} such that for $t\in I$, 
 $$
 \|(\gamma_t-\varphi_2^{-1}\alpha_t\varphi_2)|B_2\|_{cb}<2^{-3}\epsilon.
 $$
Here we again assume that $\varphi_2(B_2)\subset A_3$ as above. 
Since $B_2\supset A_2\supset \varphi_1(B_1)$, we have that
 $$
 \|(\alpha_t-\varphi_1^{-1}\varphi_2^{-1}\alpha_t\varphi_2\varphi_1)|B_1\|_{cb}
 <(2^{-2}+2^{-3})\epsilon.
 $$
By \ref{B4} we have a unitary $u_2\in A_2$ such that
 \BE
 \|(\Ad\,u_2\varphi_2\varphi_1-\id)|A_1\|&<&2^{-1}\epsilon,\\
 \max_{t\in I}\|\alpha_t(u_2)-u_2\|&<&\epsilon. 
 \EE
We choose $B_3$ for $A_3$, $\alpha$, and $2^{-3}\epsilon$ as in 
\ref{B4}. We may further assume that $B_3\ni u_2$. We then choose 
$\varphi_3\in\Aut(A)$ such that for $t\in I$,
 $$
 \|(\alpha_t-\varphi_3^{-1}\gamma_t\varphi_3)|B_3\|_{cb}<2^{-4}\epsilon.
 $$
Here again we assume that $\varphi_3(B_3)\subset A_4$. Since 
$\varphi_2(B_2)\subset A_3\subset B_3$, we have that for $t\in I$,
 $$
 \|(\gamma_t-\varphi_2^{-1}\varphi_3^{-1}\gamma_t\varphi_3\varphi_2)|B_2\|_{cb}
 <(2^{-3}+2^{-4})\epsilon.
 $$      
We then obtain a unitary $u_3\in A$ such that
 \BE
 \|(\Ad\,u_3\varphi_3\varphi_2-\id)|A_2\|&<&2^{-2}\epsilon,\\
 \max_{t\in I}\|\gamma_t(u_3)-u_3\|&<&2^{-1}\epsilon.
 \EE 
We repeat this procedure.

With $(\varphi_n)$, $(B_n)$, and $(u_n)$ obtained as above we 
proceed as follows. Let $u_1=1$, $v_1=1=v_2$, and $v_3=u_2$. We 
define a unitary $v_n$ for $n\geq 4$ by 
 $$
 v_n=u_{n-1}\varphi_{n-1}(v_{n-1}^*).
 $$
Since $v_3\in B_3$ and $\varphi_{n-1}(B_{n-1})\subset B_n$, if 
$v_{n-1}\in B_n$, we have that $v_n\in B_n$ for all $n$. By 
letting $\phi_n=\Ad\,u_n\varphi_n\Ad\,v_n^*$, we obtain an almost 
commutative diagram (cf. \cite{Ell}): 
 $$ \begin{array}{cccccc} 
 A_1& \stackrel{\rm id}{\lra}&A_3& \stackrel{\rm 
 id}{\lra}&A_5&\lra\cdots\\
 \phi_1\downarrow& 
 \phi_2\nearrow&\phi_3\downarrow&\phi_4\nearrow&\downarrow&\cdots\\
 A_2& \stackrel{\rm id}{\lra}&A_4&\stackrel{\rm 
 id}{\lra}&A_6&\lra\cdots
 \end{array}
 $$
This is almost commutative because
 \BE
 \phi_n\phi_{n-1}&=&\Ad\,u_n\varphi_n\Ad(\varphi_{n-1}(v_{n-1})u_{n-1}^*)
                     \Ad\,u_{n-1}\varphi_{n-1}\Ad\,v_{n-1}^*\\
                 &=&\Ad\,u_n\varphi_n\varphi_{n-1},
 \EE
which is in the neighborhood of the inclusion $A_{n-1}\subset 
A_{n+1}$ of norm less than $2^{-n+1}\epsilon$. Since 
$(\phi_{2n+1}(x))$ is a Cauchy sequence for $x\in \cup_nA_n$, we 
obtain a homomorphism $\phi$ of $A$ into $A$ as the extension of 
the limit of $(\phi_{2n+1})$. Since the homomorphism defined as 
the limit of $(\phi_{2n})$ is the inverse of $\phi$, it follows 
that $\phi$ is an isomorphism of $A$ onto $A$. 

For $t\in I$ we compute:
 \BE
 &&\|(\gamma_t\phi_{2n+1}-\phi_{2n+1}\alpha_t)|B_{2n+1}\| \\
 &\leq& 2\|\gamma_t(u_{2n+1})-u_{2n+1}\|
 +\|(\gamma_t\varphi_{2n+1}\Ad\,v_{2n+1}^* 
 -\varphi_{2n+1}\Ad\,v_{2n+1}^*\alpha_t)|B_{2n+1}\|\\ 
 &\leq& 2\|\gamma_t(u_{2n+1})-u_{2n+1}\|+2\|\alpha_t(v_{2n+1}^*)-v_{2n+1}^*\|
 +\|(\gamma_t\varphi_{2n+1}-\varphi_{2n+1}\alpha_t)|B_{2n+1}\| \\
  &<& 2^{-2n+2}\epsilon+2^{-2n-2}\epsilon+
 2\|\alpha_t(v_{2n+1})-v_{2n+1}\| .
 \EE
Since, for $t\in I$,
 \BE
 \|\alpha_t(v_{2n+1})-v_{2n+1}\|&\leq& \|\alpha_t(u_{2n})-u_{2n}\|
 +\|(\alpha_t\varphi_{2n}-\varphi_{2n}\gamma_t)|B_{2n}\|
 +\|\gamma_t(v_{2n}^*)-v_{2n}^*\|\\
 &\leq&2^{-2n+2}\epsilon+2^{-2n-1}\epsilon+\|\gamma_t(v_{2n})-v_{2n}\|,
 \EE                                                               
and $\|\gamma_t(v_{2n})-v_{2n}\| 
<2^{-2n+3}\eps+2^{-2n}\eps+\|\alpha_t(v_{2n-1})-v_{2n-1}\|$ and 
since $v_2=1$,  we have that for $t\in I$ 
 $$
 \|\alpha_t(v_{2n+1})-v_{2n+1}\|< 2\epsilon+2^{-2}\epsilon<3\epsilon.
 $$
From the above computations we have that for $t\in I$,
 $$
 \|(\gamma_t\phi_{2n+1}-\phi_{2n+1}\alpha_t)|A_{2n+1}\|<7\epsilon.
 $$
Thus it follows that  for $t\in I$ and $x\in \cup_nA_n$,
 $$
 \|\gamma_t\phi(x)-\phi\alpha_t(x)\|\leq 7\epsilon\|x\|.
 $$
This implies that $\alpha$ is almost conjugate to $\gamma$ and 
hence $\alpha\stackrel{cc}{\sim}\gamma$. 
\end{pf}  
  
When $B_1$ and $B_2$ are $C^*$-subalgebras of $A$ and 
$\epsilon>0$, we write $B_1\stackrel{\epsilon}{\subset} B_2$ if 
for any $x_1\in B_1$ there is an $x_2\in B_2$ with 
$\|x_1-x_2\|\leq \epsilon \|x_1\|$. We denote by $\dist(B_1,B_2)$ 
the infimum of $\epsilon>0$ such that 
$B_1\stackrel{\epsilon}{\subset} B_2$ and  
$B_2\stackrel{\epsilon}{\subset} B_1$.

Let $A$ be a UHF \cstar\ and $\alpha$ a flow on $A$. We say that 
$\alpha$ is a {\em quasi-UHF} flow if there exists a sequence 
$(A_n)$ of finite type I subfactors of $A$ such that 
$A=\ol{\cup_nA_n}$ and 
 $$
 \sup_{|t|\leq 1}\dist (A_n,\alpha_t(A_n))
 $$
converges to zero as $n\ra\infty$. Then from 6.5 of \cite{Ch} it 
follows that
 $$
 \sup_{|t|\leq1}\dist(A_n\otimes A_n,\alpha_t(A_n)\otimes\alpha_t(A_n))
 $$
converges to zero as $n\ra\infty$.

The following is an attempt to prove that the quasi-UHF flows are 
UHF. 

\begin{cor}\label{B6}
Let $A$ be a UHF \cstar\ such that $A\otimes A\cong A$ and 
$\gamma$ a universal UHF flow on $A$. If $\alpha$ is an 
approximately inner, quasi-UHF flow on $A$, then 
$\alpha\otimes\gamma$ is a cocycle perturbation of a (universal) 
UHF flow. \end{cor} 
\begin{pf}
There exists a sequence $(\eps_n)$ with $\eps_n\searrow0$ 
satisfying: For any $x\in A_n\otimes A_n$ with $\|x\|\leq1$ and 
$t\in I=[-1,1]$ there is an $x_t\in A_n\otimes A_n$ with 
$\|x_t\|\leq1$ such that 
$\|\alpha_t\otimes\alpha_t(x)-x_t\|\leq\eps_n$. 

Let $u_n$ be a self-adjoint unitary in $A_n\otimes A_n$ which 
implements the flip $\sigma|A_n\otimes A_n$.

Since, for $x\in A_n$ with $\|x\|\leq1$,
 $$ 
 (\alpha_t\otimes\alpha_t)(u_n)x\approx(\alpha_t\otimes\alpha_t)(u_nx_{-t})
 =(\alpha_t\otimes\alpha_t)(\sigma(x_{-t})u_n)\approx\sigma(x)(\alpha_t\otimes\alpha_t)(u_n),
 $$
it follows that 
$\|[(\alpha_t\otimes\alpha_t)(u_n)u_n,x]\|\leq2\eps_n\|x\|$ for 
$x\in A_n\otimes A_n$. Since there is a $v_{nt}\in A_n\otimes A_n$ 
such that $\|(\alpha_t\otimes\alpha_t)(u_n)-v_{nt}\|\leq\eps_n$, 
$\|v_{nt}\|\leq1$, $v_{nt}^*=v_{nt}$, and $\sigma(v_{nt})=v_{nt}$, 
we have that $v_{nt}u_n\in A_n\otimes A_n$ is self-adjoint and 
$\|[v_{nt}u_n,x]\|\leq4\eps_n\|x\|$ for $x\in A_n\otimes A_n$. 
Since $\|(v_{nt}u_n)^2-1\|=\|v_{nt}^2-1\|\leq2\eps_n$, $v_{nt}u_n$ 
must be close to 1 or $-1$. As we may choose $v_{nt},\ t\in I$ to 
be continuous in $t$ with $v_{n0}=u_n$, we have that $v_{nt}u_n$ 
is close to 1, i.e., $\|v_{nt}u_n-1\|\leq 6\eps_n$. Hence we get 
that $\|(\alpha_t\otimes\alpha_t)(u_n)-u_n\|\leq7\eps_n$ for $t\in 
I$. 

Since there is a sequence $(w_n)$ of unitaries in $A\otimes A$ 
such that $\sigma=\lim\Ad\,w_n$ and 
$\gamma_t\otimes\gamma_t(w_n)=w_n$, we have that the sequence 
$(u_n\otimes w_n)$ satisfies that $\sigma=\lim\Ad(u_n\otimes w_n)$ 
and 
$\|(\alpha_t\otimes\alpha_t\otimes\gamma_t\otimes\gamma_t)(u_n\otimes 
w_n)- u_n\otimes w_n\|\ra0$ uniformly in $t\in I$, i.e., the flip 
is $\alpha\otimes\gamma$-invariantly approximately inner. By the 
previous theorem we can conclude that $\alpha\otimes\gamma$ is a 
cocycle perturbation of $\gamma$. 
\end{pf}

\section{The CAR algebra}
\setcounter{theo}{0}  When $A$ is the CAR algebra (i.e., $A\cong 
M_{2^{\infty}}$), we can give a universal flow in a simple way.
 
\begin{theo}\label{C1}
Let $(\lambda_n)$ be a sequence in $\R$ and define a UHF flow 
$\alpha$ on $A=M_{2^{\infty}}$ by 
 $$
 \alpha_t=\bigotimes_{n=1}^{\infty}\Ad \left(
 \begin{array}{cc} e^{it\lambda_n}&0\\0&1\end{array}\right).
 $$
Then $\alpha$ is a universal UHF flow on $A$ if and only if there 
is a subsequence $(\lambda_{n_k})$ such that 
$\lim_k\lambda_{n_k}=0$ and $\sum_k\lambda_{n_k}^2=\infty$. 
\end{theo} 
\begin{pf}                                                
If $\alpha$ is a universal UHF flow, then the $\tau$ bottom 
marginal spectrum of $\alpha$ must be full, i.e., 
$\Gamma_{\tau,-}(\alpha)=\R_+$, where $\tau$ is the tracial state 
of $A$ (see \cite{K5}). If there is an $\eps>0$ such that 
$\sum_{n:|\lambda_n|<\eps}\lambda_n^2<\infty$, then 
$\Gamma_{\tau,-}(\alpha)\subset [\eps,\infty)$ by 5.3 of 
\cite{K5}. Hence  there must be a subsequence $(\lambda_{n_k})$ 
such that $\lim_k\lambda_{n_k}=0$ and 
$\sum_k\lambda_{n_k}^2=\infty$.  
                               
Conversely suppose that $(\lambda_n)$ satisfies the required 
property.  Let $(\mu_n)$ be another sequence in $\R$ and denote by 
$\beta$ the UHF flow constructed as in the above statement. What 
we have shown in \cite{K5} is that 
$\alpha\otimes\beta\stackrel{cc}{\sim}\alpha$ for all such 
$\beta$. Now what we have to show is that 
$\alpha\otimes\beta\stackrel{cc}{\sim}\alpha$ for any UHF flow 
$\beta$. As in \cite{K5} this follows easily from the following 
lemma. 
\end{pf} 

\ncm{\PN}{{\cal P}_N}
\begin{lem}\label{C2}
Let $(\lambda_n)$ be a sequence in $\R$ as in the above theorem 
and let T be a finite subset of $\R$ of order $2^k$. For any 
$\eps>0$ there exist an $N\in\N$ with $N>k$, a partition $P$ of 
the power set ${\cal P}_N$ of $\{1,2,\ldots,N\}$ into $2^{N-k}$ 
sets of order $2^k$, and a map $F$ of ${\cal P}_N$ onto $T$  such 
that for any $S\in P$, $F|S$ is bijective and 
 $$
 |E(x)-E(y)-F(x)+F(y)|<\eps,\ \ x,y\in S,
 $$
where $E(x)=\sum_{i\in x}\lambda_i$.               
\end{lem}
\begin{pf}
This lemma is shown in the case $k=1$ in \cite{K5}. We will extend 
the proof there to cover the general case $k>1$, at the same time 
using the result for the case $k=1$. 

We write the elements of $T$ as $\mu_1,\mu_2,\ldots,\mu_{2^k}$ in 
the increasing order. Since only the differences between $\mu_i$'s 
matter, we assume that $\mu_{2^k}=-\mu_1$. 

For a large $N$ it is certainly not difficult to find a subset 
$\{x_1,x_2,\ldots,x_{2^k}\}$ of $\PN$ such that 
$(E(x_1),E(x_2),\ldots,E(x_{2^k}))$ is almost equal to 
$(\mu_1,\mu_2,\ldots,\mu_{2^k})+C$ for some $C\in\R$ since 
$\{E(x)\ |\ x\in \PN\}$ is densely distributed in a large 
interval. We have to cover $\PN$ with the disjoint union of such 
sets $\{x_1,x_2,\ldots,x_{2^k}\}$. Then the map $F:\PN\ra T$ is an 
obvious one. 

Let $\mu_{0,i}=\mu_i$ and
 $$
 \mu_{j,i}=\frac{1}{2}(\mu_{2^j(i-1)+1}+\mu_{2^ji}),
 \ \ i=1,2,\ldots,2^{k-j}, 
 $$
for $j=1,2,\ldots,k$. Since $\mu_{2^k}=-\mu_1$, we have that 
$\mu_{k,1}=0$. For a fixed $j$, $(\mu_{j,i})$ is increasing in $i$ 
and 
 $$
 \mu_{j,2i-1}<\mu_{j+1,i}<\mu_{j,2i}.
 $$
By slightly changing $\mu_i$'s we may assume that all those 
$\mu_{j,i}$'s belong to $\Z\delta$ for some $\delta>0$, which is 
smaller than $\eps$. Furthermore, by using the result for $k=1$, 
we may assume that all $\lambda_n=\delta$. 
 
Fix a large $N$ and define $E_n(x)$ for $x\in \PN$ and $n\in\N$ 
with $n\leq N$ by
 \BE
 E_n(x)&=&\frac{1}{2}
 (\sum_{i\leq n,\,i\in x}\lambda_i-\sum_{i\leq n,\,i\not\in x}\lambda_i) 
 \\
 &=& (\sharp\{i\in x|\ i\leq n\}-n/2)\delta.
 \EE
We say that $x$ is {\em good} if there is a subsequence 
$n_1,n_2,\ldots,n_k$ of even integers such that
 $$
 E_{n_j}(x)=\mu_{k-j,\,i_j}
 $$
for some $i_j$; in particular $E_{n_k}(x)=\mu_{i_k}$. We will make 
the smallest possible choice for each $n_j$ inductively for 
$j=1,2,\ldots,k$. Since  
$\mu_{k-j-1,2i-1}<\mu_{k-j,i}<\mu_{k-j-1,2i}$, this entails 
$i_{j+1}=2i_j-1$ or $2i_j$. 

We claim that if $x$ is good, then $x$ can be matched with other 
$2^k-1$ good elements of $\PN$ in a unique way so that the 
resulting sequence $(x_1,x_2,\ldots,x_{2^k})$ satisfies that 
$(E_{n_k}(x_1),E_{n_k}(x_2),\ldots,E_{n_k}(x_{2^k}))$ exactly 
equals $(\mu_1,\mu_2,\ldots,\mu_{2^k})$. Since 
 $$x_i\bigcap 
 \{n_k+1,n_k+2,\ldots,N\}$$ 
is independent of $i$ by the construction, it follows that 
$(E(x_1),E(x_2),\ldots,E(x_{2^k}))$ is equal to 
$(\mu_1,\mu_2,\ldots,\mu_{2^k})+C$ for some $C$, where 
$E(x_i)=E_N(x_i)$. 

If $x$ is good as above, we define $x'\in \PN$ by
 $$
 x'=x\triangle\{n_{k-1}+1,n_{k-1}+2,\ldots,n_k\},
 $$
where $\triangle$ denotes difference of sets. Then $x'$ is again 
good; $E_{n_j}(x')=\mu_{k-j,i_j}$ for $j=1,2,\ldots,k-1$ and 
$E_{n_k}(x')=\mu_{\ell}$, where $\ell=i_k+1$ if $i_k$ is odd and 
otherwise $\ell=i_k-1$. In general if $(\mu_{k-j,i_j})$ is 
monotone for $j=m,m+1,\ldots,k$, we define $x'\in\PN$ by
 $$
 x'=x\triangle\{n_m+1,n_m+2,\ldots,n_k\}.
 $$
Then we have that $x'$ is good and $E_{n_k}(x')=\mu_\ell$, where 
$\ell=2^{k-m}i_m$ if $i_k$ is odd and otherwise 
$\ell=2^{k-m}(i_m-1)+1$ (since 
$\mu_{k-m,i_m}=(\mu_{2^{k-m}(i_m-1)+1}+\mu_{2^{k-m}i_m})/2$ and 
$i_k$ is either $2^{k-m}(i_m-1)+1$ or $2^{k-m}i_m$ by the 
monotonicity). If we make the smallest possible choice of $n_j'$ 
such that $E_{n_j'}(x')=\mu_{k-j,i_j'}$, we have that $n_k'=n_k$, 
where the latter is the choice for $x$. We may call this process 
the {\em reflection} at $(m,i_m)$ or loosely at $\mu_{k-m,i_m}$. 
By applying the above process inductively we can obtain the other 
$2^k-1$  good elements from one good $x$, which form the desired 
set of order $2^k$. (To understand why we get a set of $2^k$ 
elements in this way we should visualize a binary tree of depth 
$k$ with each node labeled by $\mu_{ji}$ such that the root is 
labeled by $\mu_{k1}=0$ and if a node is labeled by $\mu_{k-j,i}$ 
with $j<k$ or more precisely is addressed by $(k-j,i)$, meaning 
that it is at a distance $k-j$ from the root and at the $i$'th 
position from the left among the nodes of distance $k-j$ from the 
root, then it has two children labeled by $\mu_{k-j-1,2i-1}$ and 
$\mu_{k-j-1,2i}$ from left to right; so the leaves are labeled by 
$\mu_1=\mu_{0,1},\mu_2,\ldots,\mu_{2^k}$ from left to right. Any 
good element $x\in\PN$ corresponds to a path from the root to a 
leave of this binary tree, which is the path determined by 
$(k-j,i_j),\ j=0,1,\ldots,k$ with $i_0=1$ and $(i_j)$ as given in 
the definition of {\em good}. By the procedure indicated above we 
get $2^k-1$ elements corresponding to other $2^k-1$ paths. Note 
that to each interior node there is a reflection to be applied; 
and there are $2^k-1$ of them.) 

In this way we obtain the family $F_N$ of such sets of order 
$2^k$. Since the procedure is canonical, $F_N$ is a disjoint 
family. Let $G_N$ be the union of elements of $F_N$. If 
$x\in\PN\setminus G_N$, then $x$ is not good and in particular 
$|E_n(x)|<\mu_{2^k}$ for all $n=1,2,\ldots,N$. 

For each $k\in\Z$ with $|k|<\mu_{2^k}/\delta\equiv K$, we have 
that 
 $$
 \frac{\sharp\{x\in \PN\setminus G_N\ |\ E(x)=k\delta\}}
 {\sharp\{x\in G_N\ |\ E(x)=k\delta\}}
 $$
converges to zero as even $N$ goes to infinity. Furthermore we 
have that 
 $$
 \frac{\sharp(\PN\setminus G_N)}
 {\min_{|k|\leq K}\sharp\{x\in G_N\ |\ E(x)=k\delta\}}
 $$
converges to zero as even $N$ goes to infinity (since the ratios 
among $\sharp\{x\in G_N\ |\ E(x)=k\delta\}$ for various $k$ 
converge to 1). 

For each $x\in \PN\setminus G_N$ we specify a subset $S_x$ of 
$G_N$ of order $2K+1$  such that 
 $$
 \{E(y)\ |\ y\in S_x\}=\{k\delta\ |\ |k|\leq K\}.
 $$
If $N$ is a sufficiently large even integer, we can specify $S_x,\ 
x\in \PN\setminus G_N$ such that $S_x$'s are mutually disjoint.

Let $\{x_1,x_2,\ldots,x_{2^k}\}$ be a subset of $\PN\setminus 
G_N$. If $E(x_1)=k\delta$, we make the following substitutions 
(simultaneously): 
 $$
 y_{k-1}\leftarrow x_1,\ y_{k-2}\leftarrow y_{k-1},\ldots, 
 y_{-K}\leftarrow y_{-K+1},\ x_1\leftarrow y_{-K},
 $$
where $y_k\in S_{x_1}$ satisfies $E(y_k)=k\delta$. After these 
substitutions we have that $E(x_1)=-K\delta=\mu_1$. By making 
suitable substitutions among $\{x_i\}\cup S_{x_i}$ for 
$i=2,\ldots,2^k$ we have that $E(x_i)=\mu_i$, i.e., 
$\{x_1,x_2,\ldots,x_{2^k}\}$ satisfies the desired property. By 
these substitutions we have introduced only an error of $\delta$ 
to some sets belonging to $F_N$, i.e., for $S\in F_N$ modified by 
these processes, the elements of $S$ will be ordered as 
$x_1,x_2,\ldots,x_{2^k}$ so that $|E(x_i)-\mu_i-C|=0,\pm\delta$ 
for all $i$ and for some $C$. This completes the proof.

\end{pf}

\medskip
\small

\end{document}